\documentclass[a4paper,10pt]{article}
\title{\bf Bohr's Theorem for Monogenic Power Series}
\author{K. G\"urlebeck, J. Morais \thanks{Institute of Mathematics / Physics,
Bauhaus-University, Weimar, Germany.}}
\date{October 6, 2007}
\newtheorem{Theorem}{Theorem}[section]
\newtheorem{Definition}{Definition}[section]
\newtheorem{Remark}{Remark}[section]
\newtheorem{Proposition}{Proposition}[section]

\newtheorem{Corollary}{Corollary}[section]
\newtheorem{Lemma}{Lemma}[section]
\usepackage{amsfonts}
\newcommand{\proof}{{\it Proof.\quad}}
\begin{document}
\maketitle
\begin{abstract}
The main goal of this paper is to generalize Bohr's phenomenon from
complex one-dimensional analysis to higher dimensions in the
framework of Quaternionic Analysis.
\end{abstract}



\section{Introduction}
In 1914 Bohr discovered that there exists a radius $r \in (0,1)$
such that if a power series of a holomorphic function converges in
the unit disk and its sum has a modulus less than 1, then for
$|z|<r$ the sum of the absolute values of its terms is again less
than 1. This radius does not depend on the function.
\begin{Theorem} {\bf{(Bohr, 1914)}}
Let $f$ be a bounded analytic function in the open unit disk, with
Taylor expansion $f(z) = \displaystyle \sum_{n=0}^{\infty} a_{n}
\hspace{0.05cm} z^{n}$ convergent in the unit disk and with modulus
less than $1$. Then $\displaystyle \sum_{n=0}^{\infty} |a_{n}|
\hspace{0.05cm} r^{n} < 1$ for $0 \leq r < \frac{1}{3}$.
\end{Theorem}
This inequality known as Bohr's inequality is true for $0 \leq r <
\frac{1}{3}$ and the constant $\frac{1}{3}$ cannot be improved.

Originally, this theorem was proved for $0 \leq r <\frac16$ but soon
improved to the sharp result. In Bohr's paper \cite{HBohr1914} his
own proof was published as well as a proof by Wiener based on
function theory methods. Later, S. Sidon gave a different proof (see
\cite{SSidon1927}).

Recently, several papers were published, generalizing Bohr's theorem
to functions of $n$ complex variables (see \cite{DineenTimoney},
\cite{HD1997}, \cite{LAizenberg}). Using the standard multi-index
notations $\underline{\alpha} :=(\alpha_1,\alpha_2,...,\alpha_n) \in
\Bbb{N}_0^n$ with $|\underline{\alpha}|=\alpha_1+...+\alpha_n$,
$z:=(z_1,...,z_n)$, $z_i\in\mathbb C$,
$z^{\underline{\alpha}}:=z_1^{\alpha_1} z_2^{\alpha_2} \ldots
z_n^{\alpha_n}$, it is shown in \cite{HD1997} that if a power series
$\sum_{\underline{\alpha}} c_{\underline{\alpha}} \hspace{0.05cm}
z^{\underline{\alpha}}$ has a modulus less than $1$ in the unit
polydisc $\{ (z_1,...,z_n): \max_{1 \leq j \leq n} \hspace{0.05cm}
|z_j|<1 \}$, then the sum of the moduli of the terms is less than
$1$ in the polydisc of radius $\frac{1}{3\sqrt{n}}$.

In \cite{GueJoao22007}, the result shows the possibility to obtain a
Bohr type theorem for monogenic functions in the ball in the
Euclidean space $\Bbb{R}^3$ with the additionally condition
$f(0)=0$. It is shown that for $r<0.047$, the inequality is
satisfied. The main purpose of this paper is to check if this
theorem can be extended to all monogenic functions with
$|f({\bf{x}})|<1$ in $B_1(0)$.

Having in mind the analogy to the one-dimensional complex function
theory we want to know if the result can be proved for a ball in the
Euclidean space and not for a polydisc. It is not the goal here to
find a sharp estimate for the most general class of functions.

\section{Preliminaries}
Let $\{ {\bf{e}}_0, {\bf{e}}_1, {\bf{e}}_2, {\bf{e}}_3\}$ be an
orthonormal basis of the Euclidean vector space $\Bbb{R}^4$. The
vector ${\bf{e}}_0$ is the scalar unit while the generalized
imaginary units ${\bf{e}}_1, {\bf{e}}_2, {\bf{e}}_3$ satisfy the
following multiplication rules
\begin{eqnarray*}
& & {\bf{e}}_i {\bf{e}}_j \hspace{0.05cm} + \hspace{0.05cm}
{\bf{e}}_j {\bf{e}}_i \hspace{0.05cm} = \hspace{0.05cm} -2
\hspace{0.05cm} \delta_{i,j} \hspace{0.05cm} {\bf{e}}_0
\hspace{0.05cm}, \hspace{0.25cm} i,j = 1,2,3
\\ & & {\bf{e}}_0 {\bf{e}}_i = {\bf{e}}_i {\bf{e}}_0 = {\bf{e}}_i
\hspace{0.05cm}, \hspace{1.45cm} i = 0,1,2,3 .
\end{eqnarray*}
This non-commutative product generates the algebra of real
quaternions denoted by $\Bbb{H}$. The real vector space $\Bbb{R}^4$
will be embedded in $\Bbb{H}$ by identifying
${\bf{a}}:=(a_0,a_1,a_2,a_3) \in \Bbb{R}^4$ with the element
\begin{eqnarray*}
\mathbf{a} = a_0 {\bf{e}}_0 + a_1 {\bf{e}}_1 + a_2 {\bf{e}}_2 + a_3
{\bf{e}}_3 \in \Bbb{H} ,
\end{eqnarray*}
where $a_i$ ($i=0,1,2,3$) are real numbers. Remark that
${\bf{e}}_0=(1,0,0,0)^T$ is the multiplicative unit element of
$\Bbb{H}$ and by identifying ${\bf{e}}_0$ with $1$, it will
therefore neglected in the following notation.

The real number $\textbf{Sc} \hspace{0.04cm} \mathbf{a} := a_0$ is
called the scalar part of $\mathbf{a}$ and
$\textbf{Vec}\hspace{0.04cm}\mathbf{a} := a_1 {\bf{e}}_1 + a_2
{\bf{e}}_2 + a_3 {\bf{e}}_3$ is the vector part of $\mathbf{a}$.
Analogously to the complex case, the conjugate of $\mathbf{a}$ is
the quaternion $\overline{\mathbf{a}} := a_0 - a_1 {\bf{e}}_1 - a_2
{\bf{e}}_2 - a_3 {\bf{e}}_3$. The norm of $\mathbf{a}$ is given by
$|\mathbf{a}| = {\left(a_0^2+a_1^2+a_2^2+a_3^2\right)}^{1/2}$ and
coincides with the corresponding Euclidean norm of $\mathbf{a}$, as
a vector in $\Bbb{R}^4$. Considering the subset
\begin{eqnarray*}
\mathcal{A} := span_{\Bbb{R}}\{1,{\bf{e}}_1,{\bf{e}}_2\}
\end{eqnarray*}
of $\Bbb{H}$, the real vector space $\Bbb{R}^3$ can be embedded in
$\mathcal{A}$ by the identification of each element
$\mathbf{x}=(x_0,x_1,x_2) \in \Bbb{R}^3$ with the reduced quaternion
\begin{eqnarray*}
\textbf{x} = x_0 + x_1 {\bf{e}}_1 + x_2 {\bf{e}}_2 \in \mathcal{A}
\hspace{0,05cm}.
\end{eqnarray*}
As a consequence, we will often use the same symbol $\textbf{x}$ to
represent a point in $\Bbb{R}^3$ as well as to represent the
corresponding reduced quaternion. Note that the set $\mathcal{A}$ is
only a real vector space but not a sub-algebra of $\Bbb{H}$.

Let us consider an open set $\Omega \subset \Bbb{R}^3$ with a
piecewise smooth boundary. An $\Bbb{H}$-valued function is a mapping
$f : \Omega \longrightarrow \Bbb{H}$ such that
\begin{eqnarray*}
f(\textbf{x}) = f_0(\textbf{x}) + f_1(\textbf{x}) {\bf{e}}_1 +
f_2(\textbf{x}) {\bf{e}}_2 + f_3(\textbf{x}) {\bf{e}}_3 ,
\end{eqnarray*}
where the coordinates $f_i$ are real-valued functions defined in
$\Omega$. For continuously real-differentiable functions $f: \Omega
\longrightarrow \Bbb{H}$, the operator
\begin{equation} \label{Cauchy-Riemannoperator}
D = \partial_{x_0} + {\bf{e}}_1 \partial_{x_1} + {\bf{e}}_2
\partial_{x_2}
\end{equation}
is called the generalized Cauchy-Riemann operator. We define the
conjugate generalized Cauchy-Riemann operator by
\begin{equation} \label{conjugateCauchy-Riemannoperator}
\overline{D} = \partial_{x_0} - {\bf{e}}_1
\partial_{x_1} - {\bf{e}}_2 \partial_{x_2} .
\end{equation}
A function $f : \Omega \subset \Bbb{R}^3 \longrightarrow \Bbb{H}$ is
called $\it{left}$ (resp. $\it{right}$) $\it{monogenic}$ in $\Omega$
if
\begin{eqnarray*}
D f = 0 ~~ {\rm in} ~~ \Omega ~~ ({\rm resp}., f D = 0 ~~ {\rm in}
~~ \Omega) .
\end{eqnarray*}
From now on we only use left monogenic functions. For simplicity, we
will call them monogenic. The generalized Cauchy-Riemann operator
(\ref{Cauchy-Riemannoperator}) and its conjugate
(\ref{conjugateCauchy-Riemannoperator}) factorize the Laplace
operator in $\Bbb{R}^3$. In fact, it holds
\begin{eqnarray*}
\Delta_3 = D \overline{D} = \overline{D} D
\end{eqnarray*}
and implies that any monogenic function is also a harmonic
function.\\

From now on, we will consider the following notations: $B:=B_1(0)$
is the unit ball in $\Bbb{R}^3$ centered at the origin, $S=\partial
B$ its boundary and $d\sigma$ is the Lebesgue measure on $S$. In
what follows, we will denote by $L_2(S;\mathbb{X};\mathbb{F})$
$($resp. $L_2(B;\mathbb{X};\mathbb{F})$$)$ the $\mathbb{F}$-linear
Hilbert space of square integrable functions on $S$ $($resp. $B$$)$
with values in $\mathbb{X}$ $($ $\mathbb{X}=\mathbb{R}$ or
$\mathcal{A}$ or $\mathbb{H}$$)$, where $\Bbb{F} = \Bbb{H}$ or
$\Bbb{R}$. For any $f, g \in L_2(S;\mathcal{A};\mathbb{R})$ the
real-valued inner product is given by
\begin{eqnarray} \label{realinnerproduct}
\left< f,g \right>_{L_2(S)} = \int_{S} \textbf{Sc} (\overline{f}g)
d\sigma .
\end{eqnarray}

Each homogeneous harmonic polynomial $P_n$ of order $n$ can be
written in spherical coordinates as
\begin{eqnarray} \label{HHPsphericalcoordinates}
P_n (x) = r^n P_n(\omega), ~ \omega \in S,
\end{eqnarray}
its restriction, $P_n(\omega)$, to the boundary of the unit ball is
called $spherical$ $harmonic$ of degree $n$. From
(\ref{HHPsphericalcoordinates}), it is clear that a homogeneous
polynomial is determined by its restriction to $S$. Denoting by
$\mathcal{H}_n(S)$ the space of real-valued spherical harmonics of
degree $n$ in $S$, it is well-known (see \cite{ABR} and
\cite{CMuller}) that
\begin{eqnarray*}
\dim \mathcal{H}_n(S) = 2n + 1 .
\end{eqnarray*}
It is also known (see \cite{ABR} and \cite{CMuller}) that if $n\neq
m$, the spaces $\mathcal{H}_n(S)$ and $\mathcal{H}_m(S)$ are
orthogonal in $L_2(S;\mathbb{R};\mathbb{R})$.

Homogeneous monogenic polynomial of degree $n$ will be denoted in
general by $H_n$. In an analogously way to the spherical harmonics,
the restriction of $H_n$ to the boundary of the unit ball is called
$spherical$ $monogenic$ of degree $n$. We denote by
$\mathcal{M}_n(\Bbb{H};\Bbb{F})$ the subspace of
$L_2(B;\Bbb{H};\Bbb{F})\cap \ker D(B)$ of all homogeneous monogenic
polynomials of degree $n$. Sudbery proved in $\cite{Sud79}$ that the
dimension of $\mathcal{M}_n(\Bbb{H};\Bbb{H})$ is $n+1$. In
$\cite{DissCacao}$, it is proved that the dimension of
$\mathcal{M}_n(\Bbb{H};\Bbb{R})$ is $4n+4$.

Consider, for each $n \in \Bbb{N}_0$, a basis $\{ H^{\nu}_{n}: \nu =
1,..., \dim \mathcal{M}_n(\Bbb{H};\Bbb{F})\}$ of
$\mathcal{M}_n(\Bbb{H};\Bbb{F})$, $\Bbb{F} = \Bbb{H}$ or $\Bbb{F} =
\Bbb{R}$. Taking into account that the coordinates of $H^{\nu}_{n}$
are harmonic, for arbitrary $n, k = 0, 1,...$, we have
\begin{eqnarray} \label{relationBall/Sphere}
\left< H^{\nu}_{n}, H^{\mu}_{k} \right>_{L_2(B;\Bbb{H};\Bbb{F})} =
\delta_{n,k} \frac{1}{n+k+3} \left< H^{\nu}_{n}, H^{\mu}_{k}
\right>_{L_2(S;\Bbb{H};\Bbb{F})} .
\end{eqnarray}

\section{Homogeneous Monogenic Polynomials}
Based on the Fueter variables $\textbf{z}_1=x_1 - {\bf{e}}_1x_0$ and
$\textbf{z}_2=x_2 - {\bf{e}}_2x_0$, several systems of homogeneous
monogenic polynomials are constructed and used for different
purposes (see, e.g., \cite{BDS,Fue2,Del1070,Gue1982,HM1990,Sud79}).
Following \cite{HM1990}, being
$\underline{\gamma}=(\gamma_1,\gamma_2)$ a multi-index with
$\gamma_1+\gamma_2=n$, the generalized powers (or also Fueter
polynomials) of degree $n$ are defined by
\begin{eqnarray*}
\textbf{z}_1^{\gamma_1} \times \textbf{z}_2^{\gamma_2} &=&
\underbrace{\textbf{z}_1 \times \textbf{z}_1 \times\cdots \times
\textbf{z}_1 }_{\gamma_1\;\;times }\times{\underbrace{ \textbf{z}_2
\times \textbf{z}_2 \times \cdots \times
\textbf{z}_2}_{\gamma_2\;\;times}} \\
&=& \frac{1}{n!} \sum_{\pi(i_1,\ldots,i_{n})} \textbf{z}_{i_1}\cdots
\textbf{z}_{i_{n}}\,,
\end{eqnarray*}
\noindent where the  sum is taken over all permutations
$\pi(i_1,\ldots,i_{n})$ of
$(\underbrace{1,\cdots,1}_{\gamma_1},\underbrace{2,\cdots,2}_{\gamma_2})$.\\\\
The general form of the Taylor series of a monogenic function $f :
\Omega \subset \mathbb{R}^3 \longrightarrow \mathbb{H}$ in the
neighborhood of the origin (see, e.g., \cite{BDS,HM1990}) is given
by
\begin{eqnarray} \label{TaylorExpansion}
f = \sum_{n=0}^{\infty} \sum_{|\underline{\gamma}|=n} \left(
\textbf{z}^{\gamma_1}_1 \times \textbf{z}^{\gamma_2}_2 \right)
c_{\underline{\gamma}} ,
\end{eqnarray}
where $c_{\underline{\gamma}}=\left.\frac{1}{\gamma_1! \gamma_2!}
\hspace{0.05cm}
\partial_{x_1}^{\gamma_1} \partial_{x_2}^{\gamma_2}
f(\textbf{x})\right|_{\textbf{x}=0} \in \mathbb{H}$
are the Taylor coefficients.\\

In order to prove a collection of inequalities related to Bohr's
inequality, we need also the Fourier expansion of monogenic
functions.

In (\cite{DissCacao} and \cite{IGS2006}) $\mathbb{R}$-linear and
$\mathbb{H}$-linear complete orthonormal systems of
$\mathbb{H}$-valued homogeneous monogenic polynomials in the unit
ball of $\Bbb{R}^3$ are constructed. The main idea of these
constructions is based on the factorization of the Laplace operator.
We take a system of real-valued homogeneous harmonic polynomials and
apply the $\overline{D}$ operator to get systems of
$\mathbb{H}$-valued homogeneous monogenic polynomials. To be
precise, we introduce the spherical coordinates,
\begin{eqnarray*}
x_0 = r \cos \theta, ~ x_1 = r \sin \theta \cos \varphi, ~ x_2 = r
\sin \theta \sin \varphi,
\end{eqnarray*} where
$0 < r < \infty$, $0 < \theta \leq \pi$, $0 < \varphi \leq 2\pi$.
Each point $\textbf{x}=(x_0,x_1,x_2) \in \mathbb{R}^3 \backslash
\{0\}$ admits a unique representation $\textbf{x}=r \textbf{w}$,
where for each $i=0,1,2$ $~w_i = \frac{x_i}{r}$ and
$|\textbf{w}|=1$. Now, we apply for each $n \in \mathbb{N}_0$, the
operator $\frac{1}{2}\overline{D}$ to the homogeneous harmonic
polynomials,
\begin{eqnarray} \label{HHP}
\{ r^{n+1} U^0_{n+1}, r^{n+1} U^m_{n+1}, r^{n+1} V^m_{n+1},
m=1,...,n+1  \}_{n \in \mathbb{N}_0}
\end{eqnarray}
formed by the extensions in the ball of the spherical harmonics
(considered, e.g., in \cite{San1959}),
\begin{eqnarray} \label{sphericalharmonics}
U^0_{n+1}(\theta,\varphi) &=& P_{n+1}(\cos \theta) \nonumber \\
U^m_{n+1}(\theta,\varphi) &=& P^m_{n+1}(\cos \theta)
\cos m \varphi \\
V^m_{n+1}(\theta,\varphi) &=& P^m_{n+1}(\cos \theta) \sin m \varphi,
m=1,...,n+1 \nonumber . \cr
\end{eqnarray}
Here, $P_{n+1}$ stands for the Legendre polynomial of degree $n+1$,
given by
\begin{eqnarray*}
\left\{\begin{array}{c} P_{n+1}(t) = \displaystyle
\sum_{k=0}^{\left[ \frac{n+1}{2} \right]} a_{n+1,k} \hspace{0.15cm}
t^{n+1-2k} \\ P_0(t)=1, \hspace{1.50cm} t \in (-1,1),
\end{array}\right.
\end{eqnarray*}
with
\begin{eqnarray*}
a_{n+1,k} = (-1)^k \frac{1}{2^{n+1}} \frac{(2n+2-2k)!}{k! (n+1-k)!
(n+1-2k)!},
\end{eqnarray*}
where $[s]$ denotes the integer part of $s \in \mathbb{R}$. Also, we
stipulate this sum to be zero whenever the upper index is less then
the lower one.

The functions $P^m_{n+1}$ are called the associated Legendre
functions,
\begin{eqnarray} \label{LegendreFunctions}
P^m_{n+1}(t) := (1-t^2)^{m/2} \frac{d^m}{dt^m}P_{n+1}(t),~
m=1,...,n+1 .
\end{eqnarray}
We remark that if $m=0$, the corresponding associated Legendre
function $P^0_{n+1}$ coincides with the Legendre polynomial
$P_{n+1}$.

Notice that the Legendre polynomials together with the associated
Legendre functions satisfy several recurrence formulae. We point out
only some of them, which will be used several times in the next
section. Following \cite{LCAnd1998}, Legendre polynomials and
associated Legendre functions are solutions of a second order
differential equation, called {\it{Legendre differential equation}},
given by
\begin{eqnarray*}
(1 - t^2)(P^m_{n+1}(t))'' - 2 t (P^m_{n+1}(t))' + \left( (n + 1)(n +
2) - m^2 \frac{1}{1-t^2} \right) P^m_{n+1}(t) = 0 ,
\end{eqnarray*}
$m = 0,...,n+1$. They also satisfy the recurrence formula
\begin{eqnarray} \label{recurrenceformula}
(1 - t^2) (P^m_{n+1}(t))' = (n+m+1) P^m_n(t) - (n+1) t P^m_{n+1}(t),
\end{eqnarray}
$m=0,...,n+1$. An additional and useful identity is given by
\begin{eqnarray} \label{recurrenceformulaidentity}
P^m_m(t) = (2m - 1)!! (1 - t^2)^{m/2} ,
\end{eqnarray}
$m=1,...,n+1$.

These functions are mutually orthogonal in $L_2([-1,1])$,
\begin{eqnarray*}
\int_{-1}^1 \hspace{0.05cm} P^m_{n+1}(t) \hspace{0.05cm}
P^m_{k+1}(t) \hspace{0.05cm} dt \hspace{0.05cm} = \hspace{0.05cm} 0
\hspace{0.05cm}, n \neq k
\end{eqnarray*}
and their norms are
\begin{eqnarray*}
\int_{-1}^1 (P^m_{n+1}(t))^2 dt = \frac{2}{2n+3}
\frac{(n+1+m)!}{(n+1-m)!} , m=0,...,n+1 .
\end{eqnarray*}\

For a detailed study of Legendre polynomials and associated Legendre
functions we refer, for example, \cite{LCAnd1998} and
\cite{San1959}.\\\\ Restricting the functions of the set (\ref{HHP})
to the sphere, we obtain the spherical monogenics
\begin{eqnarray} \label{sphericalmonogenics}
X^0_n &:=& \left.\left( \frac{1}{2}\overline{D} \right) (r^{n+1}
U^0_{n+1})\right|_{r=1} \nonumber \\
X^m_n &:=& \left.\left( \frac{1}{2}\overline{D} \right) (r^{n+1}
U^m_{n+1})\right|_{r=1} \\
Y^m_n &:=& \left.\left( \frac{1}{2}\overline{D} \right) (r^{n+1}
V^m_{n+1})\right|_{r=1}, ~ m=1,...,n+1 \nonumber .
\end{eqnarray}
For each $n \in \mathbb{N}_0$, taking the monogenic extensions of
the spherical monogenics into the ball, we obtain the set of
homogeneous monogenic polynomials
\begin{eqnarray} \label{HMP}
\{r^n X^0_n, ~ r^n X^m_n, ~ r^n Y^m_n : m=1,...,n+1 \} .
\end{eqnarray}\

We need norm estimates of our functions in terms of its Taylor and
Fourier expansion are needed. In this way, we begin now to write the
homogeneous monogenic polynomials in Cartesian coordinates. In
parts, these results were already obtained in \cite{GueJoao2006} and
\cite{GueJoao2007}, without proof.
\begin{Lemma} \label{cartesiancoordinatesr^nX^l_n}
The homogeneous monogenic polynomials $r^n X^l_n (l=0,1,...,n+1)$ in
terms of Cartesian coordinates can be written as:
\begin{eqnarray*}
r^n\,X_n^l({\bf{x}}) = [r^n\,X_n^l({\bf{x}})]_0 +
[r^n\,X_n^l({\bf{x}})]_1\,{\bf{e}}_1 +
[r^n\,X_n^l({\bf{x}})]_2\,{\bf{e}}_2 ,
\end{eqnarray*}
where
\begin{eqnarray*}
[r^n\,X_n^l({\bf{x}})]_0 &=&
\sum_{k=0}^{[\frac{n-l}{2}]}\beta_{n+1,l,k}(n+1-2k-l)\,x_0^{n-2k-l}
r^{2k} \sum_{j=0}^{[\frac{l}{2}]}(-1)^j {l \choose {2j}}x_1^{l-2j}x_2^{2j}\\
                   &+ &\sum_{k=1}^{[\frac{n+1-l}{2}]}
                   \beta_{n+1,l,k}\,(2k)\,x_0^{n+2-2k-l}
                        r^{2(k-1)}
                        \sum_{j=0}^{[\frac{l}{2}]}(-1)^j
                        {l \choose 2j}x_1^{l-2j}x_2^{2j}\cr
[r^n\,X_n^l({\bf{x}})]_1 &=&
\sum_{k=1}^{[\frac{n+1-l}{2}]}\beta_{n+1,l,k}(2k)\,x_0^{n+1-2k-l}\,
r^{2(k-1)}\sum_{j=0}^{[\frac{l}{2}]}(-1)^{j+1}
{l \choose 2j}x_1^{l-2j+1}x_2^{2j} \\
& + &
\sum_{k=0}^{[\frac{n+1-l}{2}]}\beta_{n+1,l,k}\,\,x_0^{n+1-2k-l}
                         r^{2k}\sum_{j=0}^{[\frac{l-1}{2}]}(-1)^{j+1}
                         {l \choose 2j}(l-2j)x_1^{l-2j-1}x_2^{2j}\cr
[r^n\,X_n^l({\bf{x}})]_2 &=&
\sum_{k=1}^{[\frac{n+1-l}{2}]}\beta_{n+1,l,k}(2k)\,x_0^{n+1-2k-l}\,
r^{2(k-1)} \sum_{j=0}^{[\frac{l}{2}]}(-1)^{j+1}
{l \choose 2j}x_1^{l-2j}x_2^{2j+1} \\
& + & \sum_{k=0}^{[\frac{n+1-l}{2}]}
\beta_{n+1,l,k}\,\,x_0^{n+1-2k-l}
                         r^{2k} \sum_{j=1}^{[\frac{l}{2}]}(-1)^{j+1}
                         {l \choose 2j} (2j) x_1^{l-2j}x_2^{2j-1}\,,
\end{eqnarray*}
being
\begin{eqnarray*} \label{betacoefficients}
\beta_{n+1,l,k}=(-1)^k \frac{1}{2^{n+2}}\,{2n+2-2k \choose
n+1-k}{n+1-k \choose k} \,(n+1-2k)_{l-1}
\end{eqnarray*}
and $(n+1-2k)_{l-1}$ stands for the Pochhammer symbol.
\end{Lemma}
\noindent\proof Let us consider the spherical monogenics given by
(\ref{sphericalmonogenics}), explicitly described in
(\ref{sphericalharmonics}). By the definition of the Legendre
polynomials we have
\begin{eqnarray*}
P_{n+1}^{(1)}(t) = \frac{d}{dt} \sum_{k=0}^{\left[ \frac{n+1}{2}
\right]} a_{n+1,k} ~ t^{n+1-2k} = \sum_{k=0}^{\left[
\frac{(n+1)-1}{2} \right]} a_{n+1,k} (n+1-2k) ~ t^{(n+1-2k)-1} .
\end{eqnarray*}
Now, derivating recursively in order to $t$ $(l-1)$ times,
\begin{eqnarray*}
\partial_t^l ~ P_{n+1}(t) &=& P_{n+1}^{(l)}(t) \\
&=& \sum_{k=0}^{\left[ \frac{(n+1)-l}{2} \right]} a_{n+1,k} (n+1-2k)
(n+1-2k-1) \cdots (n+1-2k-(l-1)) ~ t^{(n+1-2k)-l}.
\end{eqnarray*}
By simplicity, we set
\begin{eqnarray*}
\beta_{n+1,l,k} = 2 (n+1-2k) (n+1-2k-1) \cdots (n+1-2k-(l-1)),
\end{eqnarray*}
so that, finally we get for $(\ref{LegendreFunctions})$ the
expression
\begin{eqnarray*}
P_{n+1}^{l}(\cos\theta) & =&
\sum_{k=0}^{\left[\frac{n+1-l}{2}\right]} 2 \hspace{0.05cm}
\beta_{n+1,l,k} \hspace{0.05cm} (sin\theta)^l \hspace{0.05cm}
(\cos\theta)^{n+1-2k-l} \hspace{0.05cm}.
\end{eqnarray*}
In order to express the set $\{X^l_n: l=0,1,...,n+1 \}$ in cartesian
coordinates, we consider the coordinate's relation:
\begin{eqnarray*}
&& \cos\theta = \frac{x_0}{r} \hspace{1.58cm} \cos\varphi =
\frac{x_1}{\sqrt{x_1^2+x_2^2}} \\
&& \sin\theta = \frac{\sqrt{x_1^2+x_2^2}}{r} \hspace{0.50cm}
\sin\varphi = \frac{x_2}{\sqrt{x_1^2+x_2^2}} \hspace{0.50cm}.
\end{eqnarray*}
Now, using
\begin{eqnarray*}
\cos(m\varphi) & =& \sum_{j=0}^{\left[ \frac{l}{2} \right]} (-1)^j
\hspace{0.05cm} {l \choose 2j} \hspace{0.05cm} (\cos\varphi)^{l-2j}
\hspace{0.05cm} (\sin\varphi)^{2j}
\end{eqnarray*}
and substituting in $(\ref{sphericalmonogenics})$ we obtain
\begin{eqnarray*}
r^{n+1} U_{n+1}^l ({\bf{x}}) = 2 \sum_{k=0}^
{\left[\frac{n+1-l}{2}\right]} \sum_{j=0}^{\left[ \frac{l}{2}
\right]} \beta_{n+1,l,k} ~ x_0^{n+1-2k-l} r^{2k} (-1)^j {l \choose
2j} x_1^{l-2j} x_2^{2j}.
\end{eqnarray*}
Applying the hypercomplex derivative $(\frac{1}{2}\overline{D})$ to
this expression carries our polynomials in Cartesian coordinates,
respectively. $\rule[-1ex]{2mm}{2mm}$\\\\

Similar results holds for $r^n Y^m_n$ $m=1,...,n+1)$. Let us
consider now the following function:
\begin{Definition} Let $i, j \in \Bbb{N}_0$. The function $g_{i,j}$
is given by
\begin{eqnarray*}
g_{i,j} = \left\{\begin{array}{ccccccc} 1, {\rm ~ if ~} i {\rm ~ and
~} j {\rm ~ have} {\rm ~ the ~ same ~ parity ~} \\
0 \hspace{0.05cm}, {\rm ~ if ~} i {\rm ~ and ~} j {\rm ~ have} {\rm
~ different ~ parity ~} \end{array}\right..
\end{eqnarray*}
\end{Definition}

\begin{Proposition} \label{Taylorcoefficientsr^nX^l_n}
The Taylor coefficients of the homogeneous monogenic polynomials
$r^n X^l_n (l=0,1,...,n+1)$ are given by
\begin{eqnarray*}
{\small [a^l_{\underline{\alpha}}]_0} &=& {\small g_{l,n}
\hspace{0.09cm} g_{\alpha_1,l} \hspace{0.09cm} g_{\alpha_2,0}
\hspace{0.09cm} \beta_{n+1,l,\frac{n-l}{2}}
\hspace{0.05cm} \sum_{j=0}^{[\frac{l}{2}]} (-1)^j \left(\begin{array}{cc} l \\
2j \end{array} \right)
\left(\begin{array}{cc} \frac{n-l}{2} \\
\frac{\alpha_1-l}{2}+j \end{array}\right)}
\end{eqnarray*}
\begin{eqnarray*}
&& {\small [a^l_{\underline{\alpha}}]_1 ~ = ~ g_{l-1,n} ~
g_{l-1,\alpha_1} ~ g_{\alpha_2,0} }\\
&& {\small \left[ \sum_{p=1}^{\left[\frac{n-l+1}{2}\right]}
\beta_{n+1,l,p}(2p) \left(\begin{array}{cc} p-1 \\
\frac{l-n-1}{2}+p \end{array} \right) \sum_{j=0}^{[\frac{l}{2}]}
(-1)^{j+1} \left(\begin{array}{cc} l \\
2j \end{array} \right) \left(\begin{array}{cc} \frac{n-l-1}{2} \\
\frac{\alpha_1-l-1}{2}+j \end{array} \right) \right.} \\
&+& {\small \left. \sum_{p=0}^{\left[\frac{n-l+1}{2}\right]}
\beta_{n+1,l,p} \left(\begin{array}{cc} p \\
\frac{l-n-1}{2}+p \end{array} \right) \sum_{j=0}^{[\frac{l-1}{2}]}
(-1)^{j+1} \left(\begin{array}{cc} l \\
2j \end{array} \right) (l-2j) \left(\begin{array}{cc} \frac{n-l+1}{2} \\
\frac{\alpha_1-l+1}{2}+j \end{array} \right) \right]}
\end{eqnarray*}
\begin{eqnarray*}
&& {\small [a^l_{\underline{\alpha}}]_2 ~ = ~ g_{l-1,n} ~
g_{l,\alpha_1} ~ g_{\alpha_2,1} } \\
&& {\small \left[ \sum_{p=1}^{\left[\frac{n-l+1}{2}\right]}
\beta_{n+1,l,p}(2p) \left(\begin{array}{cc} p-1 \\
\frac{l-n-1}{2}+p \end{array} \right) \sum_{j=0}^{[\frac{l}{2}]}
(-1)^{j+1} \left(\begin{array}{cc} l \\
2j \end{array} \right) \left(\begin{array}{cc} \frac{n-l-1}{2} \\
\frac{\alpha_1-l}{2}+j \end{array} \right) \right. }\\
&+& {\small \left.\sum_{p=0}^{\left[\frac{n-l+1}{2}\right]}
\beta_{n+1,l,p} \left(\begin{array}{cc} p \\
\frac{l-n-1}{2}+p \end{array} \right) \sum_{j=1}^{[\frac{l}{2}]}
(-1)^{j+1} \left(\begin{array}{cc} l \\
2j \end{array} \right) (2j) \left(\begin{array}{cc} \frac{n-l+1}{2} \\
\frac{\alpha_1-l}{2}+j \end{array} \right) \right]}
\end{eqnarray*}
\end{Proposition}
\noindent\proof The proof follows directly from Lemma
\ref{cartesiancoordinatesr^nX^l_n} by applying the partial
derivatives $\partial^{\alpha_1}_{x_1}
\partial^{\alpha_2}_{x_2}$.$\rule[-1ex]{2mm}{2mm}$\\\\
Again, we obtain analogous formulae for the Taylor coefficients of
$r^n Y^m_n$ $(m=1,...,n+1)$.
\begin{Corollary} \label{moduloscoefficients_r^nX_n^l_and_r^nY_n^m}
Let $\underline{\gamma} = (\gamma_1,\gamma_2)$ be a multi-index with
$|{\underline\gamma}| = n$. The Taylor coefficients
$a_{\underline{\gamma}}^0$, $a_{\underline{\gamma}}^m$ and
$b_{\underline{\gamma}}^m$ of the homogeneous monogenic polynomials
$r^n\,X_n^0$, $r^n\,X_n^m$ and $r^n\,Y_n^m$ satisfy the following
inequalities:
\begin{eqnarray*}
|a_{\underline{\gamma}}^0| &\leq& \frac{1}{\underline{\gamma}!}
(n+1)! \sqrt{\frac{\pi (n+1)}{2n+3}} \\
|a_{\underline{\gamma}}^m| &\leq& \frac{1}{\underline{\gamma}!}
(n+1)! \sqrt{\frac{\pi (n+1) (n+1+m)!}{2 (2n+3) (n+1-m)!}} \\
|b_{\underline{\gamma}}^m| &\leq& \frac{1}{\underline{\gamma}!}
(n+1)! \sqrt{\frac{\pi (n+1) (n+1+m)!}{2 (2n+3) (n+1-m)!}}, \qquad
m=1,...,n+1.
\end{eqnarray*}
\end{Corollary}
\noindent\proof Let $B_{r}(\textbf{x}) \subset \mathbb{R}^3$ be a
ball of radius $r$ centered at $\textbf{x}$. From \cite{GS97} we
know the Cauchy integral formula for the ball $B_1(\textbf{x})$,
\begin{equation} \label{Cauchyintegralformula}
f(\textbf{x}) = \frac{1}{4\pi} \int_{S}
\frac{\overline{\textbf{x}-\textbf{y}}}{|\textbf{x}-\textbf{y}|^3} ~
{\bf n}(\textbf{y}) f(\textbf{y}) dS_{\textbf{y}},
\end{equation}
where $\bf n$ stands for the outward pointing normal unit vector to
$S$ at $\textbf{y}$. For simplicity we just present the proof for
the homogeneous monogenic polynomials $r^nX_n^m$ $(m=1,...,n+1)$.
Applying the Cauchy integral formula to these polynomials in the
ball $B$ and taking partial derivatives with respect to $x_1$ and
$x_2$, we get
\begin{eqnarray*}
a_{\underline{\gamma}}^{m,\ast} = \frac{1}{\underline{\gamma}!}
\left. \partial_{x_1}^{\gamma_1} \partial_{x_2}^{\gamma_2}
X_n^m(\textbf{x}) \right|_{\textbf{x}=0} ~ = ~
\frac{1}{\underline{\gamma}!} \frac{1}{4\pi} \int_{S} \left.
\partial_{x_1}^{\gamma_1} \partial_{x_2}^{\gamma_2}
\frac{\overline{\textbf{x}-\textbf{y}}}{|\textbf{x}-\textbf{y}|^3}
\right|_{\textbf{x}=0} {\bf n}(\textbf{y}) X_n^m(\textbf{y})
dS_\textbf{y}
\end{eqnarray*}
taking the modulus and applying the Schwarz inequality we finally
obtain
\begin{eqnarray*}
|a_{\underline{\gamma}}^{m,\ast}| \leq \frac{1}{\underline{\gamma}!}
(n+1)! \sqrt{\frac{\pi}{2} (n+1) \frac{(n+1+m)!}{(n+1-m)!}},
\end{eqnarray*}
where $a_{\underline{\gamma}}^{m,\ast}$ denotes the Taylor
coefficients associated to the functions $X_n^m$. The previous
inequality is based on \cite{DissCacao} where the following relation
is proved
\begin{eqnarray*}
\parallel X_n^m \parallel_{L_2(S)} ~ = ~
\parallel Y_n^m \parallel_{L_2(S)} ~
= ~ \sqrt{\frac{\pi}{2} (n+1) \frac{(n+1+m)!}{(n+1-m)!}}, ~
m=1,...,n+1
\end{eqnarray*}
and on the paper \cite{KraHabil} where it was obtained that
\begin{eqnarray*}
\left\| \left.
\partial_{x_1}^{\gamma_1} \partial_{x_2}^{\gamma_2}
\frac{\overline{\textbf{x}-\textbf{y}}}{|\textbf{x}-\textbf{y}|^3}
\right|_{\textbf{x}=0} \right\|_{L_2(S)} ~ \leq ~
\frac{(n+1)!}{|\textbf{y}|^{n+2}} .
\end{eqnarray*}
Using the relation $(\ref{relationBall/Sphere})$ we get the Taylor
coefficients associated to the homogeneous monogenic polynomials
$r^n\,X_n^m$. The case $m=0$ is proved analogously. $\hspace{0.15cm}
\rule[-1ex]{2mm}{2mm}$\\\\
Besides norm estimates we also need pointwise estimates of our basis
polynomials.
\begin{Proposition} \label{modulusHMP}
For $n \in \mathbb{N}$ the homogeneous monogenic polynomials satisfy
the following inequalities:
\begin{eqnarray*}
|r^n\,X_n^0({\bf{x}})| & \leq & r^{n} (n+1)
2^n \sqrt{ \frac{\pi(n+1)}{2n+3}} \\
|r^n\,X_n^m({\bf{x}})| & \leq & r^{n} (n+1) 2^n \sqrt{\frac{\pi}{2}
\frac{(n+1)}{(2n+3)} \frac{(n+1+m)!}{(n+1-m)!}} \\
|r^n\,Y_n^m({\bf{x}})| & \leq & r^{n} (n+1) 2^n \sqrt{\frac{\pi}{2}
\frac{(n+1)}{(2n+3)} \frac{(n+1+m)!}{(n+1-m)!}}, \qquad m=1,...,n+1.
\end{eqnarray*}
\end{Proposition}
\noindent\proof Again, we prove only the case of the polynomials
$r^n\,X_n^m$ $(m=1,...,n+1)$, the proof for $r^n\,Y_n^m$ being
similar. We write these polynomials as a Taylor expansion
$(\ref{TaylorExpansion})$
\begin{eqnarray*}
r^n\,X_n^m(\textbf{x}) = \sum_{|\underline{\gamma}|=n} \left(
\textbf{z}^{\gamma_1}_1 \times \textbf{z}^{\gamma_2}_2 \right)
a^m_{\underline{\gamma}} .
\end{eqnarray*}
Consequently, modulus of $r^n\,X_n^m$ leads to
\begin{eqnarray*}
|r^n\,X_n^m(\textbf{x})| ~ \leq ~ r^{n} (n+1)! \sqrt{\frac{\pi}{2}
\frac{(n+1)}{(2n+3)} \frac{(n+1+m)!}{(n+1-m)!}} \frac{2^n}{n!}.
\end{eqnarray*}
Having in mind \cite{HM1990} we have
\begin{eqnarray*}
|\textbf{z}^{\gamma_1}_1 \times \textbf{z}^{\gamma_2}_2| \leq r^{n}
\end{eqnarray*}
for every multi-index $\underline{\gamma} = (\gamma_1,\gamma_2)$
with $|{\underline\gamma}| = n$.\\\\

For future use in this paper we will need estimates for the real
part of the spherical monogenics described in
$(\ref{sphericalmonogenics})$.
\begin{Theorem}
Given a fixed $n \in \Bbb{N}_0$, the spherical harmonics
\begin{eqnarray*}
\left\{{\bf{Sc}}\{X_{n}^0\}, {\bf{Sc}}\{X_{n}^m\},
{\bf{Sc}}\{Y_{n}^m\}: ~ m = 1,...,n \right\}
\end{eqnarray*}
are orthogonal to each other with respect to the inner product
(\ref{realinnerproduct}).
\end{Theorem}

\begin{Proposition} \label{modulusRealpartSphericalmonogenics}
Given a fixed $n \in \Bbb{N}_0$ , the moduli of the spherical
harmonics $\textbf{Sc}(X^0_n)$, $\textbf{Sc}(X^m_n)$ and
$\textbf{Sc}(Y^m_n)$ satisfy the following inequalities
\begin{eqnarray*}
|{\bf{Sc}}\{X_{n}^{l}\}| &\leq& \frac{1}{2} \frac{(n+1+l)!}{n!},
\hspace{0.37cm} l = 0,...,n
\\
|{\bf{Sc}}\{Y_{n}^{m}\}| &\leq& \frac{1}{2} \frac{(n+1+m)!}{n!},
\hspace{0.15cm} m = 1,...,n .
\end{eqnarray*}
\end{Proposition}
\noindent \proof According to the results from \cite{DissCacao}, the
real parts of the spherical monogenics are given by
\begin{eqnarray*}
\textbf{Sc}\{X_{n}^{0}\} &=& A^{0,n}(\theta) \\
\textbf{Sc}\{X_{n}^{m}\} &=& A^{m,n}(\theta) \cos(m\varphi) \\
\textbf{Sc}\{Y_{n}^{m}\} &=& A^{m,n}(\theta) \sin(m\varphi),
\end{eqnarray*}
where
\begin{eqnarray*}
A^{l,n}(\theta) = \frac{1}{2} \left( sin^{2}\theta
\frac{d}{dt}\left[P^l_{n+1}(t) \right]_{t=cos\theta} + (n+1)
cos\theta P^l_{n+1}(cos\theta) \right), ~ l=0,...,n .
\end{eqnarray*}
For simplicity sake we only present the proof for the spherical
harmonics $\textbf{Sc}(X^m_n)$ $(m=1,...,n+1)$. Making the change of
variable $t = cos\theta$ and using the recurrence formula
$(\ref{recurrenceformula})$, it follows that
\begin{eqnarray*}
\textbf{Sc}\{X_{n}^{m}\} = \frac{1}{2} (n+1+m) P^m_{n}(t).
\end{eqnarray*}
Applying the modulus in the previous expression and using the
inequality proved in \cite{Seeley1966}
\begin{eqnarray*}
|P^m_n(t)| \leq \frac{(n+m)!}{n!},
\end{eqnarray*}
for $-1 \leq t \leq 1$ and $n \geq m$, we finally obtain the
estimate
\begin{eqnarray*}
|\textbf{Sc}\{X_{n}^{m}\}| \leq \frac{1}{2} \frac{(n+1+m)!}{n!}. ~
\rule[-1ex]{2mm}{2mm}
\end{eqnarray*}\
\newline
Some of the basis polynomials described in $(\ref{HMP})$ play a
special role. Applying results from (\cite{DissCacao}, Proposition
$3.4.3$) we get:
\begin{Proposition} For $n \in \Bbb{N}_0$, the spherical monogenics
$X^{n+1}_n$ and $Y^{n+1}_n$ are given by
\begin{eqnarray} \label{monogenicconstants}
X^{n+1}_n &=& - C^{n+1,n} cos n\varphi {\bf{e}}_1 + C^{n+1,n} sin
n\varphi {\bf{e}}_2 \\
Y^{n+1}_n &=& - C^{n+1,n} sin n\varphi \hspace{0.05cm} {\bf{e}}_1 -
C^{n+1,n} cos n\varphi {\bf{e}}_2 \nonumber
\end{eqnarray}
where
\begin{eqnarray*}
C^{n+1,n} = \frac{n+1}{2} \frac{1}{sin\theta}
P^{n+1}_{n+1}(cos\theta),
\end{eqnarray*}
and their monogenic extensions into the ball belong to $ker
\overline{D}(B) \cap kerD(B)$.\\
\end{Proposition}
\begin{Remark} \label{ruleofmonogenicconstants}
The spherical monogenics $X^{n+1}_n$ and $Y^{n+1}_n$ are monogenic
constants, i.e, monogenic functions which depend only on $x_1$ and
$x_2$. Moreover, they play the role of constants with respect to the
hypercomplex differentiation $(\frac{1}{2}\overline{D})$.
\end{Remark}
\begin{Proposition} \label{modulusRealpartSphericalmonogenicsmultiplywithi}
Given a fixed $n \in \Bbb{N}_0$ , the spherical harmonics
${\bf{Sc}}(X^{n+1}_n {\bf{e}}_1)$ and ${\bf{Sc}}(Y^{n+1}_n
{\bf{e}}_1)$ are orthogonal to each other with respect to the inner
product (\ref{realinnerproduct}) and their moduli satisfy the
following inequalities
\begin{eqnarray*}
|{\bf{Sc}}\{X_{n}^{n+1} {\bf{e}}_1\}| &\leq& \frac{1}{2}
\frac{(n+1) (2n+1)!}{2^n n!} \\
|{\bf{Sc}}\{Y_{n}^{n+1} {\bf{e}}_1\}| &\leq& \frac{1}{2} \frac{(n+1)
(2n+1)!}{2^n n!}.
\end{eqnarray*}
\end{Proposition}
\noindent\proof Again, we present the proof for the spherical
harmonics $\textbf{Sc}\{X_{n}^{n+1} {\bf{e}}_1\}$, the one for
$\textbf{Sc}\{Y_{n}^{n+1} {\bf{e}}_1\}$ being similar. According to
$(\ref{monogenicconstants})$, the real part of the spherical
harmonic $X^{n+1}_n \hspace{0.05cm} {\bf{e}}_1$ is given by
\begin{eqnarray*}
\textbf{Sc}\{X_{n}^{n+1} {\bf{e}}_1\} = C^{n+1,n} \cos n\varphi.
\end{eqnarray*}
Making the change of variable $t = cos\theta$ and applying the
modulus in the previous expression, we get
\begin{eqnarray*}
|\textbf{Sc}\{X_{n}^{n+1} {\bf{e}}_1\}| &=& \frac{n+1}{2} \left|
\frac{1}{\sqrt{1-t^2}} P^{n+1}_{n+1}(t) \right|,
\end{eqnarray*}
and due to the recurrence formula
$(\ref{recurrenceformulaidentity})$ we finally obtain
\begin{eqnarray*}
|\textbf{Sc}\{X_{n}^{n+1} {\bf{e}}_1\}| = \frac{n+1}{2} \left|
\frac{1}{\sqrt{1-t^2}} (2n+1)!! (1-t^2)^{\frac{n+1}{2}} \right| \leq
\frac{1}{2} (n+1) (2n+1)!!. \hspace{0.15cm} \rule[-1ex]{2mm}{2mm}
\end{eqnarray*}\

\begin{Proposition} \label{normRealpartSphericalmonogenics}
Given a fixed $n \in \Bbb{N}_0$ , the norms of the spherical
harmonics ${\bf{Sc}}(X^0_n)$, ${\bf{Sc}}(X^m_n)$ and
${\bf{Sc}}(Y^m_n)$ are given by
\begin{eqnarray*}
\| {\bf{Sc}}(X^0_n)\|_{L_2(S)} = (n+1) \sqrt{\frac{\pi}{2n+1}}
\end{eqnarray*}
and
\begin{eqnarray*}
\| {\bf{Sc}} (X^m_n)\|_{L_2(S)} = \| {\bf{Sc}}(Y^m_n)\|_{L_2(S)} =
\sqrt{\frac{\pi}{2} \frac{(n+1+m)}{(2n+1)} \frac{(n+1+m)!}{(n-m)!}},
~ m = 1,...,n.
\end{eqnarray*}
\end{Proposition}
\begin{Proposition} \label{normRealpartSphericalmonogenicsmultiplywithi}
Given a fixed $n \in \Bbb{N}_0$ , the spherical harmonics
${\bf{Sc}}(X^{n+1}_n {\bf{e}}_1)$ and ${\bf{Sc}}(Y^{n+1}_n
{\bf{e}}_1)$ are orthogonal to each other with respect to the inner
product (\ref{realinnerproduct}) and their norms are given by
\begin{eqnarray*}
\| {\bf{Sc}}(X^{n+1}_n {\bf{e}}_1)\|_{L_2(S)} = \|
{\bf{Sc}}(Y^{n+1}_n {\bf{e}}_1)\|_{L_2(S)} = \frac{1}{2} \sqrt{\pi
(n+1) (2n+2)!}.
\end{eqnarray*}
\end{Proposition}\

\section{Bohr's Theorem}
We will denote by $X_{n}^{0,\ast},...$ the normalized basis
functions in $L_2(S;\Bbb{H};\Bbb{H})$.
\begin{Theorem} (see \cite{DissCacao}) \label{orthonormalBasis}
Let $M_n(\mathbb{R}^3; \mathcal{A})$ be the space of
$\mathcal{A}$-valued homogeneous monogenic polynomials of degree $n$
in $\Bbb{R}^3$. For each $n$, the set of $2n+3$ homogeneous
monogenic polynomials
\begin{eqnarray} \label{systemorthonormalBasis}
\left\{\hspace{0.05cm} \sqrt{2n+3} r^n\,X_{n}^{0,\ast}, \sqrt{2n+3}
r^n\,X_{n}^{m,\ast}, \sqrt{2n+3} r^n\,Y_{n}^{m,\ast}, ~ m =
1,...,n+1 \right\}
\end{eqnarray}
forms an orthonormal basis in $M_n(\mathbb{R}^3; \mathcal{A})$.
\end{Theorem}\

In \cite{GueJoao22007}, a first version of a quaternionic Bohr's
theorem was considered, therein we restricted ourselves to the case
of functions with $f(0)=0$ and we obtained an estimate in terms of a
radius of $r=0.047$.

Here, we extend our result to all monogenic functions with
$|f(\textbf{x})|<1$ in $B$, estimating a value for the radius.

\begin{Theorem}
Let $f$ be a square integrable $\mathcal{A}$-valued monogenic
function with $|f({\bf{x}})|<1$ in $B$, ${\bf{Sc}}\{f\}$ be positive
and let
\begin{eqnarray*}
\sum_{n=0}^{\infty} \sqrt{2n+3} ~ r^n \left\{ X_{n}^{0,\ast}
\alpha_{n}^{0} + \sum_{m=1}^{n+1} \left[ X_{n}^{m,\ast}
\alpha_{n}^{m} + Y_{n}^{m,\ast} \beta_{n}^{m} \right] \right\}
\end{eqnarray*}
be its Fourier expansion. Then
\begin{eqnarray*}
\sum_{n=0}^{\infty} \sqrt{2n+3} ~ r^n \left|\left\{ X_{n}^{0,\ast}
\alpha_{n}^{0} + \sum_{m=1}^{n+1} \left[ X_{n}^{m,\ast}
\alpha_{n}^{m} + Y_{n}^{m,\ast} \beta_{n}^{m} \right]
\right\}\right| < 1
\end{eqnarray*}
holds in the ball of radius $r$, with $0 \leq r <0.05$.
\end{Theorem}
\noindent\proof According to Theorem \ref{orthonormalBasis}, a
monogenic  $L_2$-function $f : \Omega \subset \Bbb{R}^3
\longrightarrow \mathcal{A}$ can be written as Fourier series
\begin{eqnarray*}
f = \sum_{n=0}^{\infty} \sqrt{2n+3} ~ r^n \left\{ X_{n}^{0,\ast}
\alpha_{n}^{0} + \sum_{m=1}^{n+1} \left[ X_{n}^{m,\ast}
\alpha_{n}^{m} + Y_{n}^{m,\ast} \beta_{n}^{m} \right] \right\},
\end{eqnarray*}
where $\alpha_{n}^{0}, \alpha_{n}^{m}$ and $\beta_{n}^{m}
(m=1,...,n+1)$ are the associated Fourier coefficients. Let us
denote by $\textbf{Sc}\{f\}$ the real part of $f$. Then,
\begin{eqnarray*}
\textbf{Sc}\{f\} &=& \frac{f + \overline{f}}{2} \\
&=& \sum_{n=0}^{\infty} \sqrt{2n+3} ~ r^n \left\{ \textbf{Sc}
\{X_{n}^{0,\ast}\} \alpha_{n}^{0} + \sum_{m=1}^{n} \left[
\textbf{Sc} \{X_{n}^{m,\ast}\} \alpha_{n}^{m} + \textbf{Sc}
\{Y_{n}^{m,\ast}\} \beta_{n}^{m} \right] \right\}.
\end{eqnarray*}
Due to Remark \ref{ruleofmonogenicconstants}, we split the function
$f$ in the following way
\begin{eqnarray*}
f &=& \sqrt{3} \alpha^0_0 X_{0}^{0,\ast} + \sqrt{3} \alpha^1_0
X_{0}^{1,\ast} + \sqrt{3} \hspace{0.05cm} \beta^1_0
Y_{0}^{1,\ast} \\
&+& \sum_{n=1}^{\infty} \sqrt{2n+3} ~ r^n \left\{ X_{n}^{0,\ast}
\alpha_{n}^{0} + \sum_{m=1}^{n} \left[ X_{n}^{m,\ast} \alpha_{n}^{m}
+ Y_{n}^{m,\ast}n\beta_{n}^{m} \right] \right\} \\
&+& \sum_{n=1}^{\infty} \sqrt{2n+3} ~ r^n \left[ X_{n}^{n+1,\ast}
\alpha_{n}^{n+1} + Y_{n}^{n+1,\ast} \beta_{n}^{n+1} \right].
\end{eqnarray*}
Based in this spliting, we introduce
\begin{eqnarray*}
f_1 &=& \frac{1}{2} \sqrt{\frac{3}{\pi}} \alpha^0_0 +
\sum_{n=1}^{\infty} \sqrt{2n+3} ~ r^n \left\{ X_{n}^{0,\ast}
\alpha_{n}^{0} + \sum_{m=1}^{n} \left[ X_{n}^{m,\ast} \alpha_{n}^{m}
+ Y_{n}^{m,\ast} \beta_{n}^{m} \right] \right\} \\
f_2 &=& -\frac{1}{2}\sqrt{\frac{3}{\pi}} \alpha^1_0 {\bf{e}}_1 -
\frac{1}{2}\sqrt{\frac{3}{\pi}} \beta^1_0 {\bf{e}}_2 +
\sum_{n=1}^{\infty} \sqrt{2n+3} ~ r^n \left[ X_{n}^{n+1,\ast}
\alpha_{n}^{n+1} + Y_{n}^{n+1,\ast} \beta_{n}^{n+1} \right],
\end{eqnarray*}
so that $f=f_1+f_2$. Then, we have
\begin{eqnarray*}
f(0) ~ = ~ f_1(0) + f_2(0)
\end{eqnarray*}
where
\begin{eqnarray*}
f_1(0) &=& \frac{1}{2} \sqrt{\frac{3}{\pi}} \alpha^0_0 \\
f_2(0) &=& -\frac{1}{2} \sqrt{\frac{3}{\pi}} \alpha^1_0 {\bf{e}}_1 -
\frac{1}{2}\sqrt{\frac{3}{\pi}} \beta^1_0 {\bf{e}}_2 .
\end{eqnarray*}

Let us assume that there exists $0<\delta<1$ such that
$|f_1|<\delta$ and $|f_2|<1-\delta$. In this way, the modulus of $f$
is preserved. We start now to study the function $f_1$. The main
idea is to compare each Fourier coefficient with the coefficient
$\alpha_{0}^{0}$. In fact, multiplying both sides of the expression
\begin{eqnarray} \label{equation}
\textbf{Sc} \{\delta - f_1\} = \delta - \textbf{Sc} \{f_1\}
\end{eqnarray}
by each real part of the homogeneous monogenic polynomials described
in $(\ref{HMP})$ and integrating over the sphere, we get these
relations. For simplicity we just present the idea applied to the
coefficients of $X_{n}^{0,\ast}$, i.e, $\alpha_{n}^{0}$. Multiplying
both sides of the expression $(\ref{equation})$ by
$\textbf{Sc}\{X_{k}^{0}\}$ and integrating, we obtain
\begin{eqnarray*}
-\sqrt{2k+3} ~ \alpha_{k}^{0} = \int_{S} \textbf{Sc} \{\delta -
f_1\} \textbf{Sc}\{X_{k}^{0}\} d\sigma
\end{eqnarray*}
with $0<\delta<1$. Now, applying the modulus we obtain finally
\begin{eqnarray} \label{alphacoefficients}
|\alpha_{k}^{0}| \sqrt{2k+3} \leq 2 \sqrt{\pi}
\frac{|\textbf{Sc}\{X_{k}^{0}\}|}{\| \textbf{Sc}\{X_{k}^{0}\}
\|_{L_2(S)}^2} \left(\delta - \frac{1}{2}\sqrt{\frac{3}{\pi}}
\alpha_0^0\right).
\end{eqnarray}
In an analogous way, we can state the following results:\\
\begin{eqnarray*}
|\alpha_{k}^{p}| \sqrt{2k+3} &\leq&  2 \sqrt{\pi}
\frac{|\textbf{Sc}\{X_{k}^{p}\}|}{\| \textbf{Sc}\{X_{k}^{p}\}
\|_{L_2(S)}^2} \left(\delta- \frac{1}{2}\sqrt{\frac{3}{\pi}}
\alpha_0^0\right). \\
|\beta_{k}^{p}| \sqrt{2k+3} &\leq& 2 \sqrt{\pi}
\frac{|\textbf{Sc}\{Y_{k}^{p}\}|}{\| \textbf{Sc}\{Y_{k}^{p}\}
\|_{L_2(S)}^2} \left(\delta- \frac{1}{2}\sqrt{\frac{3}{\pi}}
\alpha_0^0\right), ~ p=1,...,k.
\end{eqnarray*}
With some calculations and, using the Propositions
\ref{normRealpartSphericalmonogenics} and
\ref{modulusRealpartSphericalmonogenics} we finally obtain
\begin{eqnarray*}
\frac{|\textbf{Sc}\{X_{k}^{0}\}|}
{\|\textbf{Sc}\{X_{k}^{0}\}\|_{L_2(S)}^2} &\leq& \frac{1}{2\pi}
\hspace{0.05cm} \frac{(2k+1)}{k+1} \\
\frac{|\textbf{Sc}\{X_{k}^{p}\}|}
{\|\textbf{Sc}\{X_{k}^{p}\}\|_{L_2(S)}^2} &=&
\frac{|\textbf{Sc}\{Y_{k}^{p}\}|}
{\|\textbf{Sc}\{Y_{k}^{p}\}\|_{L_2(S)}^2} \leq \frac{1}{\pi}
\frac{(2k+1)(k-p)!}{(k+1+p) k!}, ~ p=1,...,k.
\end{eqnarray*}\
\newline
Finally, the previous expressions can be rewritten
\begin{eqnarray*}
|\alpha_{k}^{0}| \sqrt{2k+3} &\leq& \frac{1}{\sqrt{\pi}}
\frac{(2k+1)}{k+1} \left(\delta- \frac{1}{2}
\sqrt{\frac{3}{\pi}} \alpha_0^0\right) \\
|\alpha_{k}^{p}| \sqrt{2k+3} &\leq& \frac{2}{\sqrt{\pi}}
\frac{(2k+1)(k-p)!}{(k+1+p) k!} \left(\delta- \frac{1}{2}
\sqrt{\frac{3}{\pi}} \alpha_0^0\right) \\
|\beta_{k}^{p}| \sqrt{2k+3} &\leq& \frac{2}{\sqrt{\pi}}
\frac{(2k+1)(k-p)!}{(k+1+p) k!} \left(\delta- \frac{1}{2}
\sqrt{\frac{3}{\pi}} \alpha_0^0\right).
\end{eqnarray*}
Consequently, we can state the following inequalities:
\begin{eqnarray*}
|X^{0,\ast}_k| |\alpha_{k}^{0}| \sqrt{2k+3} &\leq&
\frac{1}{\sqrt{\pi}} (2 r)^k (2k+1) \left(\delta- \frac{1}{2}
\sqrt{\frac{3}{\pi}} \alpha_0^0\right) \\
\sum_{p=1}^{k} |X^{p,\ast}_k| |\alpha_{k}^{p}| \sqrt{2k+3} &\leq&
\frac{2}{\sqrt{\pi}} (2 r)^k (2k+1) \left(\delta- \frac{1}{2}
\sqrt{\frac{3}{\pi}} \alpha_0^0\right) \\
\sum_{p=1}^{k} |Y^{p,\ast}_k| |\beta_{k}^{p}| \sqrt{2k+3} &\leq&
\frac{2}{\sqrt{\pi}} (2 r)^k (2k+1) \left(\delta- \frac{1}{2}
\sqrt{\frac{3}{\pi}} \alpha_0^0\right).
\end{eqnarray*}\
\newline
Now, using the previous inequalities we end with
\begin{eqnarray*}
|f_1| &\leq& \frac{1}{2} \sqrt{\frac{3}{\pi}} \alpha^0_0 +
\sum_{n=1}^{\infty} \sqrt{2n+3} ~ r^n \left[ |X_{n}^{0,\ast}|
|\alpha_{n}^{0}| + \sum_{m=1}^{n} \left( |X_{n}^{m,\ast}|
|\alpha_{n}^{m}| + |Y_{n}^{m,\ast}| |\beta_{n}^{m}| \right) \right] \\
&\leq& \frac{1}{2} \sqrt{\frac{3}{\pi}} \alpha^0_0 +
\frac{5}{\sqrt{\pi}} \left(\delta- \frac{1}{2} \sqrt{\frac{3}{\pi}}
\alpha_0^0\right) \sum_{n=1}^{\infty} (2 r)^n (2n+1).
\end{eqnarray*}
Thus, we have that
\begin{eqnarray*}
|f_1| \leq \delta ~ \Longrightarrow ~ \frac{5}{\sqrt{\pi}}
\sum_{n=1}^{\infty} (2 r)^n (2n+1) \leq 1,
\end{eqnarray*}
and, the last series is convergent for $r<0.05$. In the same way, we
can study the function $f_2$. Let
\begin{eqnarray*}
f_2 = \sqrt{3} \alpha^1_0 r^0 X_{0}^{1,\ast} + \sqrt{3} \beta^1_0
r^0 Y_{0}^{1,\ast} + \sum_{n=1}^{\infty} \sqrt{2n+3} ~ r^n \left[
X_{n}^{n+1,\ast} \alpha_{n}^{n+1} + Y_{n}^{n+1,\ast} \beta_{n}^{n+1}
\right].
\end{eqnarray*}
Multiplying $f_2$ in the right side by ${\bf{e}}_1$ we get
\begin{eqnarray*}
\tilde{f_2} &:=& f_2 {\bf{e}}_1 \\
&=& \sqrt{3} \alpha^1_0 (r^0 X_{0}^{1,\ast} {\bf{e}}_1) + \sqrt{3}
\beta^1_0 (r^0 Y_{0}^{1,\ast} {\bf{e}}_1) \\
&+& \sum_{n=1}^{\infty} \sqrt{2n+3} ~  r^n \left[ (X_{n}^{n+1,\ast}
{\bf{e}}_1) \alpha_{n}^{n+1} + (Y_{n}^{n+1,\ast} {\bf{e}}_1)
\beta_{n}^{n+1} \right].
\end{eqnarray*}
We want to apply the same idea previously used for $f_1$. Taking in
consideration that $f$ is an $\mathcal{A}$-valued function, we
obtain an estimate for the coefficient $\alpha_{0}^{1}$. In a
similar way, we obtain an estimate for $\beta_{0}^{1}$ if we
multiply $f_2$ at right by ${\bf{e}}_2$. This leads to the
inequalities
\begin{eqnarray*}
|\alpha_{k}^{k+1}| \sqrt{2k+3} &\leq& 2\sqrt{\frac{\pi}{3}}
\frac{|\textbf{Sc}\{X_{k}^{k+1} {\bf{e}}_1\}|}
{\|\textbf{Sc}\{X_{k}^{k+1} {\bf{e}}_1\}\|_{L_2(S)}^2}
\left((1-\delta) -
\frac{1}{2} \sqrt{\frac{3}{\pi}} \alpha_0^1\right) \\
|\beta_{k}^{k+1}| \sqrt{2k+3} &\leq& 2\sqrt{\frac{\pi}{3}}
\frac{|\textbf{Sc}\{Y_{k}^{k+1} {\bf{e}}_1\}|}
{\|\textbf{Sc}\{Y_{k}^{k+1} {\bf{e}}_1\}\|_{L_2(S)}^2}
\left((1-\delta) - \frac{1}{2} \sqrt{\frac{3}{\pi}}
\alpha_0^1\right).
\end{eqnarray*}
being
\begin{eqnarray*}
\frac{|\textbf{Sc}\{X_{k}^{k+1}
{\bf{e}}_1\}|}{\|\textbf{Sc}\{X_{k}^{k+1} {\bf{e}}_1\}\|_{L_2(S)}^2}
= \frac{|\textbf{Sc}\{Y_{k}^{k+1}
{\bf{e}}_1\}|}{\|\textbf{Sc}\{Y_{k}^{k+1} {\bf{e}}_1\}\|_{L_2(S)}^2}
\leq \frac{2}{\pi}\frac{1}{2^n (n+1)!}.
\end{eqnarray*}
Consequently, we have proved:
\begin{eqnarray*}
|X^{k+1,\ast}_k {\bf{e}}_1| |\alpha_{k}^{k+1}| \sqrt{2k+3} &\leq&
\frac{2}{\sqrt{3 \pi}} \frac{r^k}{k!} \left((1-\delta)-
\frac{1}{2} \sqrt{\frac{3}{\pi}} \alpha_0^1 \right) \\
|Y^{k+1,\ast}_k  {\bf{e}}_1| |\beta_{k}^{k+1}| \sqrt{2k+3} &\leq&
\frac{2}{\sqrt{3 \pi}} \frac{r^k}{k!} \left((1-\delta)- \frac{1}{2}
\sqrt{\frac{3}{\pi}} \alpha_0^1 \right).
\end{eqnarray*}\
\newline
With the previous inequalities we get
\begin{eqnarray*}
|\tilde{f_2}| = |f_2| \leq \frac{1}{2} \sqrt{\frac{3}{\pi}}
\alpha^1_0 + \frac{4}{\sqrt{3\pi}} \left((1-\delta) - \frac{1}{2}
\sqrt{\frac{3}{\pi}} \alpha_0^1 \right) \sum_{n=1}^{\infty}
\frac{r^n}{n!}.
\end{eqnarray*}
Finally, we end with
\begin{eqnarray*}
|f_2| \leq 1-\delta ~ \Longrightarrow ~ \frac{4}{\sqrt{3\pi}}
\sum_{n=1}^{\infty} \frac{r^n}{n!} \leq 1,
\end{eqnarray*}
and, the last series is convergent for $r<0.56$. Finally,
\begin{eqnarray*}
\sum_{n=0}^{\infty} \sqrt{2n+3} ~ r^n \left|\left\{ X_{n}^{0,\ast}
\alpha_{n}^{0} + \sum_{m=1}^{n+1} \left[ X_{n}^{m,\ast}
\alpha_{n}^{m} + Y_{n}^{m,\ast} \beta_{n}^{m} \right]
\right\}\right| < 1
\end{eqnarray*}
converges for $0 \leq r <0.05$. $\hspace{0.15cm}
\rule[-1ex]{2mm}{2mm}$

\end{document}